\documentclass[reqno]{amsart}

\usepackage{graphics,color}

\usepackage{amssymb,amsmath,color}
\usepackage{amssymb,amsmath}
\theoremstyle{definition}
\newtheorem{defn}{Definition}[section]

\newtheorem{theorem}{Theorem}[section]
\newtheorem{lemma}[theorem]{Lemma}

\newtheorem{cor}[theorem]{Corollary}

\newtheorem{remark}[theorem]{Remark}
\newtheorem{prop}[theorem]{Proposition}

\def\l{{\lambda}}

\newcommand{\BZ}{{\mathbb Z}}

\newcommand{\al}{\alpha}
\newcommand{\G}{{\Gamma}}
\newcommand{\SL}{{SL_2(\BZ)}}

\let\a\alpha      \let\d\delta
  
  \let\l\lambda   
      
\let\GL\Lambda
\def\C{\mathbb C}
\def\G{\mathbf G}

\def\E{\mathcal E}

\def\sgn{\rm sgn}

\def\GL{{GL}}

\def\Tr{{\rm Tr}}
\def\Gal{{\rm Gal}}

\def\F{{\mathbb F}}
  
\def\Z{{\mathbb Z}}

\def\Q{{\mathbb Q}}

\def\G{\Gamma}

\def\ta{\widetilde{A}}

\DeclareRobustCommand{\qed}{%
  \ifmmode 
  \else \leavevmode\unskip\penalty9999 \hbox{}\nobreak\hfill
  \fi
  \quad\hbox{$\square$} \\}

\begin{document}
\thanks{ The author was
supported in part by an NSF-AWM mentoring travel grant for women.
She would like to thank the Pennsylvania State University  for
hosting her visit during May-June 2006.}

\title{On Atkin and Swinnerton-Dyer congruence relations (3)}

\author{Ling Long}
\address{Department of Mathematics\\Iowa State University\\Ames, IA 50011 \\USA}
\begin{abstract} In the previous two papers with the same title (\cite{lly05} by W.C. Li, L. Long, Z. Yang and
\cite{all05} by A.O.L. Atkin, W.C. Li, L. Long), the authors have
studied special families of  cuspforms for noncongruence arithmetic
subgroups. It was found that the Fourier coefficients of these
modular forms at infinity satisfy three-term Atkin and
Swinnerton-Dyer congruence relations which are the $p$-adic analogue
of the three-term recursions satisfied by the coefficients of
classical Hecke eigenforms.

In this paper, we first consider   Atkin and Swinnerton-Dyer type
congruences which generalize the three-term congruences above. These
weaker congruences are satisfied by cuspforms for special
noncongruence arithmetic subgroups. Then we will exhibit an infinite
family of noncongruence cuspforms, each of which satisfies
three-term Atkin and Swinnerton-Dyer type congruences for almost
every prime $p$. Finally, we will study a particular space of
 noncongruence cuspforms. We will show that the attached
$l$-adic Scholl representation is isomorphic to the $l$-adic
representation attached to a classical automorphic form. Moreover,
for each of the four residue classes of odd primes modulo 12 there
is a basis so that the Fourier coefficients of each basis element
satisfy  three-term Atkin and Swinnerton-Dyer congruences in the
stronger original sense.
\end{abstract}

\maketitle
\section{Introduction}
This paper is a continuation of two previous  papers with the same
title: \cite{lly05} by W.C. Li, L. Long, Z. Yang and \cite{all05} by
A.O.L. Atkin, W.C. Li, L. Long. Here,  we continue to explore
modular forms for noncongruence arithmetic subgroups.  The serious
study of these functions was initiated in the late 1960's by Atkin
and Swinnerton-Dyer \cite{a-sd} and further developed by Scholl
\cite[etc]{sch85b,sch86,sch88}. In particular, under a general
assumption Scholl has established a system of compatible $l$-adic
representations of the absolute Galois group attached to each space
of noncongruence cuspforms. In \cite{lly05, all05}, two intricate
cases have been exhibited in which noncongruence cuspforms are
related closely to classical congruence automorphic forms as
follows: From the $p$-adic point of view, in each of these cases
there is a simultaneous (or semi-simultaneous) basis for almost all
primes $p$ such that each basis function is an ``eigenform" for the
``$p$-adic" Hecke operators $T_p$ (c.f. Section \ref{sec: p-adic}).
Moreover, the ``traces" of these $T_p$  are Hecke eigenvalues (up to
at most an ideal class character).  From the representation point of
view, in each case the associated $l$-adic Scholl representation
$\rho_l$ can be decomposed into a direct sum of 2-dimensional
subrepresentations, when $\rho_l$ is restricted to a suitable Galois
subgroup. Furthermore, each subrepresentation is shown to be
isomorphic to an $l$-adic representation attached to a classical
congruence automorphic form. In this paper, we intend to give a more
general discussion on the first aspect and a further discussion on
the second one.

In this paper, we continue to consider modular forms for
noncongruence character groups (following the notation of  A.O.L.
Atkin). To be precise, a  character group is an arithmetic subgroup
 $\G$ which is normal in its congruence closure $\G^0$ (the
smallest congruence subgroup which contains $\G$) with abelian
quotient. For example, the lattice groups studied by Rankin
\cite{rankin67} and all groups in \cite{lly05, all05, Kurth-Long06}
are character groups and most of them are noncongruence. The
construction of these groups endows the corresponding modular forms
with special arithmetic properties.

This paper is organized in the following way: Section \ref{sec:ncc}
is devoted to the properties of modular forms for character groups.
In Section \ref{sec:Gn}, we  study the arithmetic of weight 3
cuspforms for a special family of noncongruence character groups
denoted by $\G_n$.   In particular, we will prove the following
theorem
\begin{theorem}
  For every positive integer $n$ and almost every prime $p$, the space  $S_3(\G_n)$ of weight 3 cuspforms
  for $\G_n$
  has a rational basis independent of $p$ such that for every basis element
  $\sum_{n\ge 1} a(n) w^n$, there exist a natural number $r$
  depending on $p$ and $n$ and two algebraic numbers $A(p)$ and $B(p)$ with
  $|A(p)|\le \binom{2r}{r}p^r$ and $|B(p)|\le p^{2r}$ such that for
  any positive integer $n$
  \begin{equation}\label{eq:asdtype}
\frac{a(np^r)-A(p)a(n)+B(p)a(n/p^r)}{(pn)^2}
  \end{equation} is integral at $p$.
\end{theorem}

The congruence relation \eqref{eq:asdtype} is weaker but more
general than those Atkin and Swinnerton-Dyer relations obtained in
\cite{lly05, all05}.


In  Section \ref{sec:G6}, we investigate the arithmetic properties
of weight 3 cuspforms for $\G_6$ in the above family. This case is
very similar to the one obtained in \cite{all05}. The main result is
\begin{theorem} Let $\rho_{3,l,6}: \Gal(\overline{\Q}/\Q)
\rightarrow \text{Aut}(W_{3,l,6})$ be the $l$-adic Scholl
representations attached to weight 3 cuspforms for $\G_6$. Then
$\rho_{3,l,6}$ are isomorphic up to semisimplification to the
$l$-adic representations attached to a classical congruence
automorphic form. Moreover, for every prime $p\nmid 6l$, the space
of weight 3 cuspforms has a basis, depending on the congruence class
of $p$ modulo 12, such that for each basis element $\sum_{n\ge 1}
a(n) w^n$ and all $n\in \mathbb N$
 \begin{equation}
\frac{a(np)-A(p)a(n)+B(p)a(n/p)}{(pn)^2}
  \end{equation} is integral at some place in $\Z[i]$ above $p$ where $A(p)$ and $B(p)$ are two
algebraic numbers  with  $|A(p)|\le2p$ and $|B(p)|\le p^{2}$.
Moreover, $A(p)$ and $B(p)$ can be obtained from the coefficients of
 a congruence cuspform.
\end{theorem}

In this paper, we use $\G$ to denote a finite index subgroup of
$\SL$ such that $-I_2\notin \G$ and use $\G^0$ to denote its
congruence closure.
 Let $X_{\G}$
denote the compact modular curve for $\G$. We assume all modular
curves considered here are defined over $\Q$. Let $\frak{M}_{\G}$
denote the field of meromorphic modular functions for $\G$, i.e. the
field of rational functions on $X_{\G}$. In particular,
$\frak{M}_{\G}$ is a finite extension of $\frak{M}_{\G^0}$. Let $k$
be a positive integer. By $S_k(\G)$ we mean the space of weight $k$
holomorphic cuspforms for $\G$. In the sequel, let $\omega_n=e^{2\pi
i/n}$ and $L_n=\Q(\omega_n)$. For any field $K$, let
$G_K=\Gal(\overline{K}/K)$. Unless otherwise mentioned, we follow
all other notation used in \cite{all05}.

\section{Modular forms for  character groups}\label{sec:ncc}
\subsection{Character groups}\label{subsec:ncc}

\begin{defn}\label{defn:ncc}
 An arithmetic subgroup $\G$  is called a character group of another arithmetic subgroup $\G^0$ of if $\Gamma$ is normal in $\G^0$ with
 abelian
 quotient.
\end{defn}

Since the quotient group $\G^0/\G$ acts on the $\C$-vector space
$S_k(\G)$ by the stroke operator $|$, by the representation of
finite abelian group, $S_k(\G)$ can be decomposed into a direct sum
of representation subspaces parameterized by the characters of
$\G^0/\G$.

\subsection{Character groups with cyclic $\G^0/
\G$} By the Fundamental Theorem of Finite Abelian Groups, we will
focus on  character groups with $\G^0/\G$  cyclic
 in the sequel. Assume $\G$ is a  character group
and $\G^0/\G=\langle \zeta \G\rangle \cong \Z_n$ where
$n=[\G^0:\G]$.  By the representation theory of finite cyclic
groups,
$$S_k(\G)=\bigoplus _{1\le j \le n} S_k(\G^0, \chi^j),$$ where $\chi(\zeta \G)=\omega_n$ and $\chi^j$ ($1\le j \le n$) are the  characters of $\G^0/\G$.
Here $S_k(\G^0, \chi^j)$ consists of functions $f$ in $S_k(\G)$ such
that $f|_{\zeta}=\omega_n^j f$. Thus, when $(j,n)>1$, $f\in
S_k(\G^0,\chi^j)$ is  a modular form for an intermediate group
sitting between $\G$ and $\G^0$. So, $\bigoplus_{(j,n)>1}
S_k(\G^0,\chi^j)$ consists of ``old" forms for supergroups of $\G$;
while $\bigoplus_{(j,n)=1} S_k(\G^0,\chi^j)$ consists of ``new"
forms orthogonal to those ones in $S_k(\G)^{\text{old}}$.
Accordingly, we denote these two spaces by $S_k(\G)^{\text{old}}$
and $S_k(\G)^{\text{new}}$ respectively. Under the Petersson inner
product, we have
$$S_k(\G)=S_k(\G)^{\text{new}}\bigoplus S_k(\G)^{\text{old}}.$$

\subsection{Special symmetry} Since $\zeta $ normalizes $\G$, it induces an order
$n$ cyclic covering map from the modular curve $X_{\G}$  of $\G$ to
the modular curve $X_{\G^0}$ of $\G^0$. We use $\zeta$ again to
denote such an involution and assume its minimal field of definition
 to be $\Q(\omega_n)$. Moreover, we assume that  as a
simple finite extension of $\frak{M}_{\G^0}$, $ \frak{M}_{\G}$ is
generated by $t_n$ with minimal polynomial $(t_n)^n-t$ for some
 $t\in \frak{M}_{\G^0} $ with rational Fourier coefficients
at infinity.  Then $X_{\G}$ satisfies a symmetry
\begin{equation}\label{eq:zetamap}
   \zeta : t_n\mapsto \omega_n^{-1} t_n.
\end{equation}   Such a map $\zeta $ is
defined over $L_n=\Q(\omega_n)$.

\subsection{$l$-adic Scholl representations}Assume $\dim _{\C}
S_k(\G)=d$. For any prime number $l$, let  $\rho_{k,l}: G_{\Q}
\rightarrow \text{Aut}(W_{k,l})$ be the $l$-adic Scholl
representation attached to $S_k(\G)$ which is unramified outside of
a few primes \cite{sch85b}. Here $W_{k,l}$ is a $2d$-dimensional
$\Q_l$-vector space. Denoted by $N$ the product of all ramifying
primes of $\rho_{k,l}$ except $l$ together with all prime divisors
of $n$. Let $p\nmid Nl$ be a prime and $F_p$ be the canonical
Frobenius conjugacy class in the quotient of the decomposition group
at $p$ by the inertia group at $p$. Scholl has also shown that the
characteristic polynomial of $\rho_{k,l}(F_p)$ for any $p\nmid Nl$
has integral coefficients and the eigenvalues have the same absolute
value $p^{(k-1)/2}$.

When $\G$ is a  character group as above. The map $\zeta $
\eqref{eq:zetamap} endows $W_{k,l}\otimes \Q_l(\omega_n)$ with
 an order $n$ involution. Under our assumptions,
$W_{k,l}\otimes \Q_l(\omega_n)$ decomposes into eigenspaces of
$\zeta$. Let $\mathcal{W}_j$ be such an eigenspace with eigenvalue
$\omega_n^j$. Let
$$W_{k,j}^{\text{new}}=\bigoplus_{(j,n)=1} \mathcal{W}_j \quad \mathrm{and}
\quad W_{k,j}^{\text{old}}=\bigoplus_{(j,n)>1} \mathcal{W}_j.$$  The
space $W_{k,l}^{\text{new}}$ is simply denoted by $W^{\text{new}}$
when there is no ambiguity.

\begin{lemma}\label{lem:Fxi} Assume $p\nmid Nl$ is a prime. In that
case
  $$F_p\zeta =\zeta ^p F_p,\quad  \text{and} \quad \zeta  F_p= F_p \zeta ^{\hat{p}},$$ where $\hat{p}$
  denotes the inverse of $p$ in $(\Z/n\Z)^{\times}$.
\end{lemma}\begin{proof}Modulo the prime ideal $(p)$, we
have \begin{eqnarray*} F_p\zeta  (t_n)&=&F_p(\omega_n^{-1}
t_n)=\omega_n^{-p}t_n^p=\zeta ^p F_p(t_n).\\
\zeta  F_p(t_n)&=&\omega_n^{-1} t_n^p=(\omega_n^{-\hat{p}}
t_n)^p=F_p \zeta ^{\hat{p}}(t_n).
\end{eqnarray*}

\end{proof}

Consequently, for any $w\in \mathcal{W}_j$, $\zeta (F_p w)=
F_p\zeta^{\hat p}w=\omega_n ^{\hat p} (F_p w).$ Thus $F_p
\mathcal{W}_j\subseteq \mathcal{W}_{j\hat{p}}.$ By iteration, $F_p$
permutes $\mathcal{W}_j$ and

\begin{equation*}\label{eq:FpW}
F_p \mathcal{W}_j=\mathcal{W}_{j\hat{p}}.
\end{equation*}

 Since  the group generated by all $F_p$ with $p\nmid Nl$ acts on $\{\mathcal{W}_j\}_{j\in
(\Z/n\Z)^{\times}}$ transitively,
\begin{equation}
 \dim _{\Q_l(\omega_n)}\mathcal{W}_j ={\d_p},
\end{equation}
where $\d_p$ is an integer independent of $j$ in $
(\Z/n\Z)^{\times}$. Therefore
\begin{equation}
\dim _{\Q_l(\omega_n)}W^{\text{new}}= \d_p \cdot \phi(n),
\end{equation} where $\phi(n)$ is the
Euler number of $n$. For any $p\nmid n$, use $O_n(p)$ to denote the
order of $p$ in the multiplicative group $(\Z/n\Z)^{\times}$. Let
$r=O_n(p)$. Moreover, we have
 \begin{equation}
\mathcal L_{j,p}=\bigoplus_{m=1}^{r} \mathcal{W}_{j\hat{p}^{m}}
 \end{equation} is
 invariant under $F_p$.

\begin{cor}\label{cor:WdecompasGalmodul}
  The spaces $W_{k,l}^{\text{new}}$ and $W_{k,l}^{\text{old}}$ are  invariant
  under each $F_p$ when $p\nmid Nl$. Therefore as a
$\Gal(\overline{\Q}/\Q)$-module, $W_{k,l}\otimes
\Q_l(\omega_n)=W_{k,l}^{\text{new}}\bigoplus W_{k,l}^{\text{old}}$.
\end{cor}
\begin{proof}
For any $p\nmid Nl$, $(\hat{p},n)=1$, and $(j,n)=(j \hat{p},n)$. The
assertions follow naturally.
\end{proof} In the remaining part of this section, we will focus on $W_{k,l}^{\text{new}}$
and let $\rho_{k,l}^{\text{new}}: G_{\Q}\rightarrow
\text{Aut}(W_{k,l}^{\text{new}})$.

 As we have extended the
scalar field to include $\omega_n$, there exists a basis $\mathcal
B$ of $W^{\text{new}}$ under which the matrix of $\zeta $ on
$W^{\text{new}}$ is  a diagonal matrix of the following block form
\begin{equation}\label{eq:zeta}\zeta=
\begin{pmatrix}
    \omega_n^{j_1}I_{{\d_p}}&0&0&\cdots&0&0\\
0&\omega_n^{j_2}I_{{\d_p}}&0&\cdots&0&0\\
\cdots&\cdots&\cdots&\cdots&\cdots&\cdots\\ 0&0&0&\cdots&\omega_n^{j_{\phi(n)-1}}I_{{\d_p}} &0\\
0&0&0&\cdots&0&\omega_n^{j_{\phi(n)}}I_{{\d_p}}
  \end{pmatrix},
\end{equation} where $j_i$ runs through the set $ (\Z/n\Z)^{\times}$.

If $B$ is an operator on a vector space $L$, let
$\text{Char}(L,B)(T)$  denote the characteristic polynomial of $B$
on $L$ with variable $T$.

\begin{lemma}\label{lem:Lj}
 Let $p\nmid Nl$ be a prime number   and  $j \in (\Z/n\Z)^{\times}$. Let $r=O_n(p)$.
  Then $$\text{Char}(\mathcal L_{j,p}, F_p)(T)\in \Q_l(\omega_n)[T^r].$$
\end{lemma}
\begin{proof}
  Assume that under $\mathcal B$ the matrix of $F_p$ on $W^{\text{new}}$ is  $(E_{i,j})_{1\le i,j\le \phi(n)}$ where  $E_{i,j}$
  are
  ${\d_p}\times {\d_p}$ matrices. The commutativity $F_p\zeta =\zeta ^{p}F_p$ in Lemma \ref{lem:Fxi} and \eqref{eq:zeta} imply that
  $E_{i,j}=(0)_{{\d_p}\times {\d_p}}$ are all-zero matrices unless $i\cdot(-p)+j=0 \text{ mod }n.$ So the matrix of $(F_p)^n$ restricted to
  $\mathcal L_{j,p}$  consists of block forms with diagonal blocks  $(0)_{{\d_p}\times {\d_p}}$ unless
  $r|n$. Therefore the trace of $F_p^n$ on
  $\mathcal L_{j,p}$ is 0 unless $r|n$.  Consequently the
  characteristic polynomial of $F_p$ on $\mathcal L_{j,p}$ is  in terms of $T^r$.
\end{proof}

\begin{lemma}\label{lem:charLj}
$\text{Char}(\mathcal
L_{j,p},F_p)(T)=\text{Char}(\mathcal{W}_j,F_p^r)(T^r)$.
\end{lemma}
\begin{proof}
  By Lemma \ref{lem:Lj}, if $\al$ is a solution of
$\text{Char}(\mathcal L_{j,p},F_p)(T)=0$, so is $\omega_r \al$.
Hence,  all the roots of $\text{Char}(\mathcal L_{j,p},F_p)(T)$ are
$\{\omega_r^j\al_m\}_{1\le j\le r, 1\le m\le {\d_p}}$.

On the other hand $\text{Char}(\mathcal{W}_j,F_p^r)(T)$ coincides
with the characteristic polynomial of the ${\d_p}\times{\d_p}$
matrix $\prod_{m=1}^rE_{p^{-m}j, p^{-m+1}j}$. Since  the order of
the product does not effect the characteristic polynomial,
$$\text{Char}(\mathcal L_{j,p},F_p^r)(T)=\left (\text{Char}(\mathcal{W}_j,F_p^r)(T)\right )^r=\prod_{1\le m\le {\d_p}} (T-\al_m^r)^r.$$
Hence $\text{Char}(\mathcal{W}_j,F_p^r)(T) =\prod_{1\le m\le {\d_p}}
(T-\al_m^r)$. It follows
$$\text{Char}(\mathcal{W}_j,F_p^r)(T^r) =\prod_{1\le m\le {\d_p}}
(T^r-\al_m^r)=\text{Char}(\mathcal L_{j,p},F_p)(T).$$

\end{proof}

\begin{prop}\label{cor:WinTr}
The polynomial  $\text{Char}(W^{\text{new}},F_p)(T)$ is in $\Z[T^r]$
with all roots having
 the same
  absolute value $p^{(k-1)/2}$.
\end{prop}
\begin{proof} By Lemma \ref{lem:Lj} and  Scholl's Theorem which says $\text{Char}(W_{k,l},F_p)(T)\in
\Z[T]$ and its eigenvalues have the same absolute value
$p^{(k-1)/2}$, it suffices to show
$\text{Char}(W^{\text{new}},F_p)(T)\in \Z[T]$. Let $\G'$ be an
intermediate group $\G\subset \G'\subset \G^0$. By our assumption,
$X_{\G'}$ has a model over $\Q$. Let $W_{k,l}(\G')$ be the $l$-adic
Scholl representation space associated with $S_k(\G')$. Then
$\text{Char}(W_{k,l}(\G'),F_p)(T)\in \Z[T]$. When $m=[\G^0:\G']$ is
a prime,
$$\text{Char}(W^{\text{old}}_{k,l}(\G'),F_p)(T)=\text{Char}(W_{k,l}(\G^0),F_p)(T) \in
\Z[T]$$ and hence $$\text{Char}(W^{\text{new}}_{k,l}(\G'),F_p)(T)
\in \Z[T].$$ By induction, for any divisor $m|n$,
$\text{Char}(\bigoplus_{j,(j,n)=n/m}\mathcal W_j,F_p)(T) \in \Z[T]$.
In particular,
$$\text{Char}(W_{k,l}^{\text{new}},F_p)(T)=\text{Char}(\bigoplus_{j,(j,n)=n/n}\mathcal W_j,F_p)(T)\in
\Z[T^r].$$
\end{proof}

\subsection{Induced representations}\label{sec:inducedrepn}
As a consequence of the above discussions we have for any $p\nmid
Nl$ $$\Tr(\rho_{k,l}^{\text{new}}(F_p))=0 \quad  \text{if} \quad  p
\neq 1 \text{ mod }n.$$ Since the images of all Frobenius elements
will determine any Galois representation up to semisimplification,
this implies that the character of $\rho_{k,l}^{\text{new}}$ is
invariant under any twisting by Dirichlet characters of modulus
divisible by $n$. Let $\mathcal W_j$ as above, $L_n=\Q(\omega_n)$,
and
 $\rho_j: G_{L_n} \rightarrow
\text{Aut}(\mathcal W_j)$ for any $j\in (\Z/n\Z)^{\times}$. The
group $G_{L_n}$ is an index $\phi(n)$ subgroup of $G_{\Q}$ and
$\mathcal{W}_j$ is a ${\d_p}$-dimensional representation space of
$G_{L_n}$. Applying a result in \cite[Prop. 19]{seree-linrep}, we
obtain the following result.
\begin{prop}\label{thm:induced}Let $j\in
(\Z/n\Z)^{\times}$. We have
$$\rho_{k,l}^{\text{new}}=\mathrm{Ind}_{G_{\Q(\omega_n)}}^{G_{\Q}}\rho_j.$$
\end{prop}

\subsection{$p$-adic spaces and Atkin and Swinnerton-Dyer type congruences}\label{sec: p-adic}Let
$S_k(\G, \Z_p)= S_k(\G)\otimes \Z_p,$ and $V$ be the $p$-adic Scholl
space attached to $S_k(\G)$ with $S_k(\G,\Z_p)$ as a subspace
\cite{sch85b}. These
 $p$-adic spaces are endowed with the action of the Frobenius
morphism, which is denoted by $F$. The operator $\zeta$ also acts on
$V$ and its order is still $n$. Like before, we can define $V_j$ to
be the $\omega_j$-eigenspace of $\zeta$ on $V$ (over
$\Q_p(\omega_n)$) as well as $V^{\text{new}}$ and $V^{\text{old}}$
similarly.

In \cite{a-sd}, Atkin and Swinnerton-Dyer  observed that some
noncongruence modular forms satisfy three-term Hecke-like recursions
in a $p$-adic sense. Their observations have been verified by
Cartier for weight 2 cases \cite{car71} and Scholl \cite{sch85b} for
all 1-dimensional cases. Here, we will first define Atkin and
Swinnerton-Dyer type congruence relations.
\begin{defn}
  Let $K$ be a number field and $f=\sum_{n\ge 1} a(n)w^n$ a weight $k$ modular form. The function $f$ is
  said to satisfy the {Atkin and Swinnerton-Dyer type congruence relation} at
  $p$ given by polynomial $T^{2d}+A_1T^{2d-1}+\cdots +A_{2d}\in K[T]$ if
  for all $n\in \Z$,
  $$\frac{a(np^d)+A_1 a(n p^{d-1})+\cdots+A_da(n)+\cdots+A_{2d}
  a(n/p^d)}{(np)^{k-1}}$$ is integral at some place of $K$ above the
  prime number $p$.
\end{defn}

To examine three-term Atkin and Swinnerton-Dyer congruence relations
computationally, Atkin uses ``$p$-adic" Hecke operators which we
will explain. For any prime $p\nmid N$, define the $p$-adic Hecke
operator to be $T_p=U_p+s(p)\cdot p^{k-1}V_p$ which acts on a weight
$k$ form $ f=\sum a(n)w^n$ as follows:
$$ f|_{T_p}= \sum_n (a(np)+s(p)p^{k-1}a(n/p) )  w^n.$$  Here, $s(p)$ is an algebraic number with absolute
value 1.
 It is called a ``$p$-adic" Hecke operator as the $n$th
coefficient of $f|_{T_p}$ is only determined up to modulo some
suitable power of $p$ depending on $n$.  In general, it is totally
nontrivial to guess the right $s(p)$ values so that $T_p$ has
``eigen" forms in $S_k(\G)$.  Moreover $s(p)$ are rarely Dirichlet
characters. Unlike the classical Hecke operators which can be
diagonalized simultaneously, these operators are a priori defined
over various $p$-adic fields. It is clearly exceptional that the
space $S_k(\G)$ has a simultaneous eigen basis for all ``$p$-adic"
Hecke operators $T_p$, which is the case in \cite{lly05}.

\subsection{Atkin and Swinnerton-Dyer type congruences satisfied by modular forms for noncongruence character groups}
 Applying  Theorem 5.6 in \cite{sch85b}, we have
\begin{theorem}\label{thm:scholl1}[Scholl, \cite{sch85b}] For any prime $p\nmid Nl$.
 $$ \text{Char}(V, F)(T)=\text{Char}(W,F_p)(T)\in \Z[T].$$
\end{theorem}

\begin{prop}[Fang et al \cite{L5} ]\label{prop:F5}Let $G$ be a finite group of automorphisms of the
elliptic modular surface $\pi: \E_{\G} \rightarrow X_{\G}$
associated with $\G$. Let $\chi$ be an irreducible character of $G$
and if $V$ is a representation of $G$, let $V^{\chi}$ denote the
$\chi$-isotypical subspace of $V$. Let $K/\Q$ be  the field of
definition of the representation whose character is $\chi$ and let
$\lambda$ be a place of $K$ above $l$ and $\wp$ be a place of $K$
above $p$. Finally, we let $r$ to be the smallest positive integer
such that $(F_p)^r \in G_K$. Then
\begin{equation}
  \text{Char}((V \otimes K_{\wp})^{\chi}, F^r)(T)=\text{Char}((W \otimes K_{\l})^{\chi},
  (F_{p})^r)(T).
\end{equation}
\end{prop}

In our case, we use $G=\G^0/\G$ which is generated by $\zeta \G $.
Then $(V \otimes K_{\wp})^{\chi^j}=V_j$ and $(W \otimes
K_{\l})^{\chi^j}=\mathcal{W}_j$.

\begin{cor}
  If $j\in (\Z/n\Z)^{\times}$, then $\dim _{\Q_p(\omega_n)} V _j ={\d_p}$ where ${\d_p}$ is  even.
\end{cor}

\begin{cor}\label{cor:Wnew=Vnew}
$$\text{Char}(V^{\text{new}},F)(T)=\text{Char}(W^{\text{new}},F_p)(T)\in \Z[T].$$
Moreover, the
 roots of these monic polynomials have the same absolute value
 $p^{(k-1)/2}$.
\end{cor}

\begin{theorem}\label{thm:asdtype} Assume $\G$ is a character group of $\G^0$ with $\G^0/\G=<\zeta \G>\cong \Z_n,$ $X_{\G}$ has a model
defined over $\Q$, and the action of $\zeta$ on $X_{\G}$ is defined
over $\Q(\omega_n)$.
   Let $p\nmid Nl$ be a prime, $j\in (\Z/n\Z)^{\times}$, and
$r=O_n(p)$. Assume the space $V_j$ has a basis consisting of forms
with Fourier coefficients in $L_n$. Let $H_j(T)=\text{Char}(\mathcal
L_j,F_p)(T)\in L_n[T^r]$. For every $f\in V_j$ with coefficients in
$L_n$, $f$ satisfies an Atkin and Swinnerton-Dyer type congruence at
$p$ given by  $H_j(T)$.
\end{theorem}
\begin{proof}
It follows from  Scholl's arguments in \cite{sch85b}.
\end{proof}
 In particular, when ${\d_p}=2$, we have
\begin{equation}\label{eq:3termasdtype}
  H_j(T)=T^{2r}-A_j(p)T^r+B_j(p),
\end{equation}where $ \quad |A_j(p)|\le
\binom{2r}{r}p^{(k-1)r/2}, |B_j(p)|=p^{(k-1)r}$. The corresponding
three-term recursion is weaker than the original three-term Atkin
and Swinnerton-Dyer congruence relation when $r>1$.

\section{Atkin and Swinnerton-Dyer type congruences satisfied by cuspforms for $\G_n$}\label{sec:Gn}
In this section, we will use the results obtained in the previous
section to derive three-term Atkin and Swinnerton-Dyer type
congruence relations satisfied by weight 3 cuspforms for $\G_n$.
\subsection{The groups $\G_n$}
 In
\cite{lly05, all05}, the following family of noncongruence character
group is considered.

Let $$\G^1(5)=\{ \gamma\in \SL \big| \gamma \equiv \begin{pmatrix} 1&0\\
* &1
\end{pmatrix} \text{ mod }5 \}.$$ It is an index 12 congruence subgroup of
$\SL$ with 4 cusps $\infty, -2, 0,-5/2$. This group has four
generators $\gamma_{\infty}, \gamma_{-2}, \gamma_0,\gamma_{-5/2}$
(each stabilizes one cusp as indicated by the subscripts)  subject
to one relation
$\gamma_{\infty}\gamma_{-2}\gamma_{0}\gamma_{-5/2}=I_2$. Let
$\varphi_n$ be the homomorphism
\begin{eqnarray*}
\G^1(5)&\rightarrow & \C^{\times}\\
\gamma_{\infty}&\mapsto &\omega_n\\
\gamma_{-2}&\mapsto& \omega_n^{-1}\\
\gamma_0, \gamma_{-5/2}&\mapsto &1.
\end{eqnarray*}  The kernel $\G_n$ of $\varphi_n$ is an index $n$ noncongruence character group of
 $\G^1(5)$.  When $n\neq 1,5$,  $\G_n$ is a
noncongruence subgroup.

Let $E_1, E_2$ be the weight 3 Eisenstein series for $\G^1(5)$ as in
\cite{lly05}. A Hauptmodul for $\G^1(5)$ is $t=\frac{E_1}{E_2}$
which generates $\frak M_{\G^1(5)}$ and
$t_n=\sqrt[n]{\frac{E_1}{E_2}}$ is a Hauptmodul for $\G_n$. Such a
group $\G_n$ has two distinguished normalizers in $\SL$:
$\zeta=\begin{pmatrix} \label{eq:AonE} 1&5\\0&1
\end{pmatrix}, \quad A=\begin{pmatrix}
  2&-5\
  \\1&-2
\end{pmatrix}$. In
particular
\begin{equation}
  E_1|_{A}=E_2, \quad E_2|_{A}=-E_1, \quad  \quad t|_{A}=-\frac{1}{t}.
\end{equation}

We use the following algebraic model for the elliptic modular
surface $\E_{\G_n}$ associated to the group $\G_n$ to investigate
the roles of the above operators.
\begin{eqnarray}\label{eq:E_n}
\begin{split}
\E_{\G_n}: \quad y^2=t_n^n(x^3-\frac{1+12(t_n^{n}-t_n^{3n})+14t_n^{2n}+t_n^{4n}}{48t_n^{2n}}x\\
  +\frac{1+18(t_n^{n}-t_n^{5n})+75(t_n^{2n}+t_n^{4n})+t_n^{6n}}{864t_n^{3n}}).
 \end{split}
\end{eqnarray}The actions of $A$ and $\zeta$ on $\E_{\G_n}$ are:
\begin{equation}\label{eq:A}
A(x,y,t_n)=(-x,\frac{y}{t_n},\frac{\omega_{2n}}{t_n});
\end{equation}
\begin{equation}
  \zeta(x,y,t_n)=(x,y,\omega_{n}^{-1} t_n).
\end{equation}
Therefore, $A$ and $\zeta$ are defined over $\Q(\omega_{2n})$ and
$\Q(\omega_{n})$ and have order 4 and $n$ respectively. In fact, the
$A$ map used in \cite{lly05} which is defined over $\Q$ is a
derivation of the $A$ map here.

Following the notation used before, let
$$\rho_{k,l,n}: G_{\Q} \rightarrow \text{Aut}(W_{k,l,n})$$ be
the $l$-adic Scholl representation attached to $S_k(\G_n)$.

\subsection{The space $S_3(\G_n)$}
When $n>1$, $k=3, j\in (\Z/n\Z)^{\times},$ $$ \quad {\d_p}=\dim
_{\Q_l(\omega_n)}\mathcal W_j=\dim _{\Q_p(\omega_n)} V_j=2.$$   For
any positive integer $n$, let $h_j^{[n]}=\sqrt[n]{E_1^{n-j}E_2^j}\in
\Z[1/n][w]], w=e^{2\pi i z/5n}$. By  \cite[Prop 2.1]{all05},
$S_3(\G_n)=\langle h_j^{[n]}\rangle _{j=1}^{n-1}$. It is
straightforward to verify that $h_j^{[n]} \in V_j$.  According to
Theorem \ref{thm:asdtype}, we have
\begin{theorem}For every basis element $h_j^{[n]}$ of $S_3(\G_n)$ and any  prime $p>n$,
$h_j^{[n]}$ satisfies   a three-term Atkin and Swinnerton-Dyer type
congruence relation at $p$ given by $\text{Char}(\mathcal W_j,
F_p^r)(T^r)$ as in Lemma \ref{lem:charLj} where $r=O_n(p)$.
\end{theorem}
To achieve this result, we only need  the operator $\zeta$. When we
consider in addition the $A$ map, we will obtain the following
result using \cite[Prop. 19]{seree-linrep}.

Let  $\rho_{3,l,n}^{\text{new}}: G_{\Q} \rightarrow
\text{Aut}(W_{3,l,n}^{\text{new}})$ be the $l$-adic representation
attached to the cuspforms genuinely belonging to $\G_n$ as before.
When $n=2,3,4$, let $B$ be $A$, $A$, and $\zeta$; let $K$ be $\Q(i),
\Q(\sqrt{-3}),$ and $\Q(i)$ respectively.
   The action of $B$ on $W_{3,l,n}^{\text{new}}$ satisfies $B^2=-1$.
   Decompose $W_{3,l,n}^{\text{new}}$ according to $(\pm i)$-eigenspaces of $B$. As a Galois $G_K$-module,
   $W_{3,l,n}^{\text{new}}=W_{3,l,n,i}^{\text{new}}\oplus W_{3,l,n,-i}^{\text{new}}$. Let
   $\rho_{3,l,n,\pm i}^{\text{new}}: G_K \rightarrow
   \text{Aut}(W_{3,l,n,\pm i}^{\text{new}})$. Then we have
\begin{prop}\label{rem:serre}
  When $n=2,3,4$,
  $$\rho_{3,l,n}^{\text{new}}=\mathrm{{Ind}} _{G_{K}}^{G_{\Q}}
  \rho_{3,l,n,\pm i}^{\text{new}}.$$
\end{prop}

\begin{remark}
  The above theorem provides a different perspective for the results
  in \cite{lly05}. When $n=2$, each $\rho_{3,l,n,\pm i}^{\text{new}}$ is a
  one-dimensional representation, hence a character of $G_{\Q(i)}$. Naturally it  corresponds to
  a cuspform $\eta(4z)^6$  with complex
  multiplication (\cite[section 8]{lly05}). When $n=3$, $\rho_{3,l,3}^{\text{new}}$ is  induced  from a 2-dimensional representation of
  $G_{\Q(\sqrt{-3})}$. Such a point of view was suggested by J.P.
  Serre as one of his comments on \cite{lly05}. A similar result also
  holds for the $n=6$ case when another derivation of $A$, called $\widetilde{A}$
   \eqref{eq:twistA} is used.
\end{remark}

\section{$S_3(\G_6)$}\label{sec:G6}
When $k=3$, $\dim_{\Q_l(\omega_n)} W_{3,l,n}^{\text{new}}=2\phi(n)$.
Accordingly, only when $n=2,3,4,$ or $6$ do we have
$\dim_{\Q_l(\omega_n)} W_{3,l,n}^{\text{new}}\le 4$. Since the first
three cases have been handled
 in \cite{lly05, all05}, we will now treat the $n=6$ case here.
The current case under consideration shares a lot of similarities
with the case $n=4$ studied in \cite{all05}.  To avoid duplication,
we will refer the readers to \cite{all05} for some arguments. In
this case, we use the following action which is a variation of $A:$
\begin{eqnarray}\label{eq:twistA}
  \widetilde{A}(x,y,t_6)&=&(-x,
  \frac{y}{t_6}, \frac{i}{t_6}).
\end{eqnarray} We use
 $\widetilde{A}$ here since it is defined over a smaller
field $\Q(i)$ while $A$ is defined over $\Q(e^{2\pi i/12})$.

Associated with $S_3(\G_6)$ is a 10-dimensional $l$-adic Scholl
representation \cite{sch85b} $$\rho_{3,l,6}: \Gal(\overline{\Q}/\Q)
\rightarrow \text{Aut}(W_{3,l,6}).$$

In this case, we can pick $N=6$ (c.f. Section
\ref{sec:ladicnumerical}). As before, the space
$W_{3,l,6}\otimes_{\Q_l} \Q_l(\omega_6)$ decomposes naturally into
eigenspaces of $\zeta$. Denoted by $\mathcal{W}_j$ the
$\Q_l(\omega_6)$-eigenspace of $\zeta$ with eigenvalue $\omega_6^j$.
Similar to the discussion in Section 5.3 of \cite{all05},
$\mathcal{W}_3$ (resp. $\mathcal{W}_2\oplus \mathcal{W}_4$) is
isomorphic to the $l$-adic Scholl representation associated with
weight 3 cuspforms for $\G_2$ (resp. $\G_3$) tensoring with
$\Q_l(\omega_6)$. The remaining piece $\mathcal{W}_1\oplus
\mathcal{W}_5$, denoted by $W_{3,l,6}^{\text{new}}$ (or simply
$W^{\text{new}}$), is a representation space of
$\Gal(\overline{\Q}/\Q)$ (c.f. Cor. \ref{cor:WdecompasGalmodul}). By
Prop. \ref{cor:WinTr}, for every Frobenius element $F_p$ where
$p\nmid 6$,
$$\text{Char}(W_{3,l,6}^{\text{new}},F_p)(T)\in \Z[T].$$

 Let $V$ be the $p$-adic Scholl space attached to $S_3(\G_6)$ with the action of Frobenius $F$.
We will drop the superscript $[6]$ below for convenience. Decompose
$V$ accordingly into $V^{\text{old}}\oplus V^{\text{new}} $ such
that
$$ h_2, \
h_3 , \ h_4  \in V^{\text{old}} \text{ and } V^{\text{new}} =V_1
\oplus V_5,$$ where $ h_1 \in V_1, \ $ and $ h_5   \in V_5$.

By considering $h_j=E_2\cdot {t_n}^{n-j}$ and the actions of $\zeta$
and $\widetilde{A}$ on $E_1, E_2$, and $t_n$, we obtain
$$h_j|\zeta=\omega_n^j h_j, \quad  h_j|\widetilde{A}=-i^{n-j}h_{n-j}.$$
It follows
\begin{equation}
  \zeta  \ta\zeta =\ta.
\end{equation}
A straightforward computation reveals that
\begin{equation}\label{eq:Fta}
  F\ta F^{-1}=\zeta ^{3(p-1)/2} \ta.
\end{equation}

\subsection{Numerical data } In this section, we will indicate how
the results were discovered numerically.

\subsubsection{$p$-adic
side}\label{sec:numericalAtkin and Swinnerton-Dyer} Since $h_3$
(resp. $h_2$ and $h_4$) has been discussed in \cite{lly05} (resp.
\cite{all05}), we will investigate  three-term Atkin and
Swinnerton-Dyer congruences for $h_1$ and $h_5$ only. Using the
``$p$-adic" Hecke operators mentioned in Section \ref{sec: p-adic},
A.O.L. Atkin observed that for any odd prime $p>3$,
\begin{itemize}
  \item $h_1|_{T_p}= c\cdot h_1, \quad h_5|_{T_p}=c_2\cdot h_5,$ when  $p=1 \,(\text{mod
  }6)$
with eigenvalues $c$  and $c_2$.

  \item $h_1|_{T_p}= c\cdot h_5,\quad  h_5|_{T_p}=c_2\cdot h_1,$  when $  p=5\, (\text{mod }6)$, where   $c=c_2$ if
$p=1 \text{ mod }4$ (resp. $c=-c_2$ if $p=-1 \text{ mod }4$). Hence
$h_1\pm h_5$ (resp. $h_1\pm i h_5$) are eigenfunctions with
eigenvalues $\pm \sqrt{c\cdot c_2}$.

\end{itemize}

Here,  $s(p)=1$, if $p=1 \text{ mod }12$ and $s(p)=-1$ otherwise.

  In table 1, we list  constants $c$ canonically modified as follows:
  \begin{itemize}
\item when $p=1 \text{ mod }12$,  we divide  $c$ by $(-3)^{(p-1)/4} \text{ mod }p=\pm
1$. \end{itemize}

 \noindent   Otherwise we divide  $c$ by a
canonical square root $c_1$ determined as follows:

\begin{itemize}
\item  when $p=5 \text{ mod }12$,  $c_1=\pm\sqrt{-1}=- (-3)^{(p-1)/4} \text{ mod }p;$
\item when $p=7 \text{ mod }12$,  $c_1=\pm\sqrt{-3}= (-3)^{(p+1)/4} \text{ mod }p;$
\item when $p=11 \text{ mod }12$, $c_1=\pm \sqrt{3}= (-3)^{(p+1)/4} \text{ mod }p.$

  \end{itemize}

\begin{center}
$$\begin{tabular}{|c|c|c||c|c|c|}
\hline
$p$& $p \text{ mod }12$& modified $c$ &$p$& $p \text{ mod }12$& modified $c$\\
\hline
5 &5& 7&37 &1& $-10$\\
7 &7& 5&41 &5 &$-50$ \\
11& 11& $-5$&43 &7 &$-10$\\
13& 1& 20&47 &11 &$-50$\\
17 &5 &$-8$&53 &5& $-47$\\
19 &7& $-6$&59 &11 &20 \\\
23 &11& 2&61 &1 &$-64$\\
29 &5& 10&67 &7& $-50$\\
31 &7& 31&71 &11 &0
\\
\hline
\end{tabular}
$$
(Table 1)
\end{center}

 \subsection{$l$-adic representation
side}\label{sec:ladicnumerical} Like the discussions in
\cite[Section 5]{lly05} and \cite[Section 3.1]{all05}, there is an
$l$-adic representation $\rho_{3,l,6}^*$ constructed from the
elliptic modular surface $\E_{\G_6}$. Like \cite[Section
3.1]{all05}, we can show that $\rho_{3,l,6}^*$ is unramified at 5
and $\rho_{3,l,6}$ is isomorphic to $\rho_{3,l,6}^*$. The
computation on $\rho_{3,l,6}^*$ is very explicit. By using a Magma
program, we obtain the following list:
\begin{center}
$$\begin{tabular}{|c|cl|}
\hline $p$&$\text{Char}(W^{\text{new}},F_p)(T)$&\\ \hline
5&$(T^2+7iT-5^2)(T^2-7iT-5^2)$&$=T^4-T^2+5^4$\\
7&$(T^2+5\sqrt{-3}T-7^2)(T^2-5\sqrt{-3}T-7^2)$&$=T^4-23T^2+7^4$\\
11&$(T^2-5\sqrt{-3}T-11^2)(T^2+5\sqrt{-3}T-11^2)$&$=T^4-167T^2+11^4$\\
13&$(T^2+20T+13^2)^2$&\\
17&
$(T^2+8iT-17^2)(T^2-8iT-17^2)$&$=T^4-514T^2+17^4$\\
19&$(T^2+6\sqrt{-3}T-19^2)(T^2-6\sqrt{-3}T-19^2)$&$=T^4-614T^2+19^4$\\
23&$(T^2-2\sqrt{-3}T-23^2)(T^2+2\sqrt{-3}T-23^2)$&$=T^4-1046T^2+23^4$\\
37&$(T^2+10T+37^2)^2$&\\
 53&$(T^2-47iT-53^2)(T^2+47iT-53^2)$&$=T^4-3409T^2+53^4$\\\hline
\end{tabular}
$$
(Table 2)
\end{center}

\subsection{Congruence side}\label{subsec:g} We have identified a
level 108 weight 3 congruence newform by using Magma with the first
few terms given as below
\begin{eqnarray*} \widetilde{g}&:=& q + (\frac{1}{10}u - \frac{17}{10})q^2 + (-\frac{1}{5}u - \frac{3}{5})q^4 +
7q^5 + (-\frac{1}{2}u + \frac{7}{2})q^7\\&& + 8q^8  + (\frac{7}{10}u
- \frac{119}{10})q^{10} + (\frac{1}{2}u - \frac{7}{2})q^{11} +
20q^{13}\\&& + (\frac{1}{2}u + \frac{23}{2})q^{14}  + (\frac{4}{5}u
- \frac{68}{5})q^{16} - 8q^{17} + (\frac{3}{5}u -
\frac{21}{5})q^{19}\\&& + (-\frac{7}{5}u - \frac{21}{5})q^{20} +
(-\frac{1}{2}u - \frac{23}{2})q^{22}
\\&&+ (-\frac{1}{5}u + \frac{7}{5})q^{23}+\cdots-10q^{37}+\cdots- 47q^{53}+\cdots
\\&=&\sum _{n\ge 1} a(n) q^n,
\end{eqnarray*}
where $u$ is a root of  $x^2 - 14x + 349=0$. Different choices of
$u$ give rise to two weight 3 congruence newforms which are the same
up to a twisting by $\chi_{-4}$. We have verified that up to $p=541$
the $p$th coefficients listed in Table 2 and beyond match the
corresponding modified coefficients of the congruence newform (with
$u=7-10 \sqrt{-3}$). In particular, in the comparison, we use
$a(p)/(i\sqrt{3})$ instead of $a(p)$ when $p=3 \text{ mod }4$.

\section{Modularity of $\rho_{3,l,6}$}\label{sec:modularity} The reader
will find the following section resembles   Section 8 of
\cite{all05}. But the modularity is achieved via a method similar to
that used in \cite{lly05} since Livn\'e's criterion does not apply
to the current case. Moreover, we will focus on the modularity of
$\rho_{3,l,6}^{\text{new}}$ in this section.

Let $\l=1+i$.   Denote by $\rho_{\widetilde{g}}'$ the $l$-adic
Deligne representation attached to the newform $\widetilde{g}$
above. Let $L=\Q(i,\sqrt[4]{3})$. Like \cite[Section 8.1]{all05},
define an ideal class character $\chi$ by composing the Artin
reciprocity map  from the idele class group of $\Q(i)$ to
$\Gal(L/\Q(i))$ with the isomorphism which sends $\Gal(L/\Q(i))$ to
$\langle i\rangle $. In particular,
\begin{itemize}
\item $\chi(p)=1$ for any $p=3 \text{ mod }4$;
\item  $\chi(p)=\pm 1$ for any prime $v$ above $p$ when
$p=1 \text{ mod }12$, depending on the value of $(-3)^{(p-1)/4} \mod
v$;
\item $\chi(p)=\pm i$ for any prime $v$ above $p$ when $p=5 \text{ mod }12$,
depending on the value of $(-3)^{(p-1)/4} \text{ mod }v$.
\end{itemize}

Denote by $\rho_{\widetilde{g}}$ the restriction of
$\rho_{\widetilde{g}}'$ to  $G_{\Q(i)}$. There is a cusp form $g$
for $GL_2(\Q(i))$ corresponding to $\rho_{\widetilde{g}}$. The form
$g$ is the lifting of $\widetilde{g}$  over $\Q(i)$ under the base
change by Langlands. Since $\widetilde{g}$ and $\widetilde{g}$
twisted by $\chi_{-4}$ both lift to  $g$, corresponding to
$\rho_{\widetilde{g}} \otimes \chi$ is the cusp form $g_{\chi}$,
called $g$ twisted by $\chi$, for $GL_2(\Q(i))$. The local
$p$-factors of the L-function $L(s,g_{\chi})$ of $g_{\chi}$ are very
similar to those for $f_{\chi}$ in \cite[section 8]{all05}. The
function $L(s,g_{\chi})$ does not dependent on the two possible
choices for $\chi$. Moreover, the local $p$-factors  of
$L(s,g_{\chi})$ agree with $\text{Char}(W^{\text{new}}, F_p)(T)$
when $T$ is replaced by $p^{-s}$ for all primes $p$ listed in Table
2.

The action of $\widetilde{A}$  is defined over $\Q(i)$. When
restricted to $\Q(i)$, $\rho_{3,l,6}^{\text{new}}: G_{\Q}
\rightarrow \text{Aut}(W^{\text{new}})$ decomposes  into two
representations $\rho^{\text{new}}_{\pm}: G_{\Q(i)} \rightarrow
\text{Aut}(W^{\text{new}}_{\pm })$ according to the $(\pm
i)$-eigenspaces of $\widetilde{A}$. Each $W^{\text{new}}_{\pm}$ is a
2-dimensional $\Q_{1+i}$-vector space.

\begin{theorem}\label{thm:Ca}  The representation $W^{\text{new}}_+$ (or $W^{\text{new}}_-$) of $G_{\Q(i)}$,  is isomorphic, up to semisimplification, to the
representation space for $\rho_{\widetilde{g}} \otimes \chi$.
\end{theorem}

Unlike the case in \cite{all05}, we cannot apply Livn\'e's criterion
as the traces are not always even. Like the discussion in
\cite{all05}, we use a special case of the Faltings-Serre's
criterion \cite{ser_1}. The ramified places are $1+i$ and $3$.

\begin{theorem}[Serre]\label{thm:serre}
Let $\rho_1$ and $\rho_2$ be representations of
$\Gal(\overline\Q/\Q(i))$ to $\GL_2(\Z[i]_{1+i})$. Assume they
satisfy the following two conditions:
\begin{itemize}
\item[(1)] $\det(\rho_1)=\det(\rho_2)$; \item[(2)] the two
homomorphisms from $\Gal(\overline\Q/\Q(i))$ to
    $\GL_2(\F_2)$, obtained from the reductions of $\rho_1$ and
    $\rho_2$ modulo $1+i$, are surjective and equal.
\end{itemize}
If $\rho_1$ and $\rho_2$ are not isomorphic, then there exists a
pair $(\widetilde G, t)$, where $\widetilde G$ is a quotient of the
Galois group $\Gal(\overline\Q/\Q(i))$ isomorphic to either
$S_4\times \{\pm 1\}$ or $S_4$ or $S_3\times \{\pm 1\}$, and the map
$t: \widetilde G \to \F_2$ has value 0 on the elements of
$\widetilde G$ of order $\le 3$, and 1 on the other elements.
\end{theorem}

The idea of Faltings-Serre's criterion has been recast in
\cite{lly05}.

Since $\l=1+i$, $\Z_{2}(i)/\l\Z_{2}(i)\cong \F_2$.
\begin{lemma}\label{lem:det}
  As 1-dimensional representations of $G_{\Q(i)}$, $\det
  \rho^{\text{new}}_{+}\cong \det \rho_{\widetilde{g}} \otimes \chi.$
\end{lemma}
\begin{proof}
We use the idea of  the proof of  \cite[Lemma 5.2]{lly05}. The only
quadratic extensions of $\Q(i)$ unramified outside of $1+i$ and $3$
are $\Q(i)(\sqrt{d})$ where $d=i,(1+i),(1-i),3, 3i, 3(1+i),3(1-i)$.
Respectively, $p=3+2i,3+2i,3-2i,4+i,3+2i,4-i,4+i$ are inert in
$\Q(i)(\sqrt{d})$. Hence, it suffices to compare the values of the
two characters $\det
  \rho^{\text{new}}_{+}$ and $\rho_{\widetilde{g}} \otimes \chi$ at primes $3\pm 2i$
  and $4\pm i$.
\end{proof}

\begin{lemma}
There is only one representation $\rho$ from  $G_{\Q(i)} $
 to $\GL_2(\F_2)$ unramified at $1+i$ and $3$ such that the characteristic
 polynomials of $\rho(F_{\wp})$ are equal to those given in Table 2 modulo
 $\lambda$. Furthermore $\rho$ is surjective.

\end{lemma}
\begin{proof}
   Let ${\sgn}:  SL_2(\F_2) \cong S_3 \rightarrow \{ \pm1 \}$
   be the sign function.
  The kernel of $\rho\circ \sgn$ fixes a quadratic field $\Q(i)(\sqrt{d})$ of $\Q(i)$
  which is unramified outside of $1+i$ and $3$. The possible $d$ are
$1,i, 1+i, 1-i, 3,3i,3+3i,3-3i.$ By the assumption, the image of
$F_{7\pm 2i}$ under $\rho$ is $T^2+T+1$. Hence its image under
$\rho$ is of order $3$ and under $\rho\circ \sgn$ its image is 1.
However, when $d=i,1+i,1-i,3,3+3i,3-3i$, either $F_{7+2i}$ or
$F_{7-2i}$ is inert in $\Q(i)(\sqrt{d})$. So the only possible
choices are $d=1,3i$.

 Now we will show $\rho$ is surjective. If  $d=1$, the  kernel of $\rho\circ \sgn$ fixes a
  cyclic cubic Galois extension $M$ of $\Q(i)$ which is unramified outside of
  $1+i$ and 3. Hence $\Gal(M/\Q)=\Z_6$ or $S_3$. Since the order of $\rho(F_7)$ is at least 6, hence  $\Gal(M/\Q)=\Z_6$ where $M=\Q(i)(\omega_9+\omega_9^{-1})$
  is the splitting field of $f(x)=x^3-3x+1$ over $\Q(i)$. The polynomial $f(x)$ is irreducible  modulo
  $3+2i$. Meanwhile, the characteristic polynomial of
  $F_{3\pm 2i}$ on $W^{\text{new}}$ modulo $\l$ is $T^2+1$. Its order should
  be 1 or 2. This leads to a contradiction.

 Finally assume that $\ker(\rho\circ \sgn)=\Gal(\overline{\Q}/K)$ for some $S_3$ extension of $\Q(i)$ which contains $\Q(i)(\sqrt{3i})$. I.e.   $K$ is the  splitting field
of an irreducible cubic polynomial $f(x)=x^3+px+q,\quad p,q,\in
\Z[i]$ such that a) the discriminant $-4p^3-27q^2=\pm
(3i)(1+i)^{2\al}3^{2\beta}, \a,\beta\in \Z_{\ge 0};$ b) $f(x)$ is
reducible  modulo $3+2i$, $6+i$ (their characteristic polynomials in
the residue field are $T^2+1$ and hence of order at most 2) and is
irreducible when modulo $7+2i$ (its characteristic polynomial in the
residue field is $T^2+T+1$ and hence of order 3). Such a field can
be determined uniquely and is the splitting field of
$f(x)=x^3-3x-2i$. So  $\rho$ is unique up to isomorphism.

\end{proof}

  It is easy to see $\rho^{\text{new}}_+$ and $\rho_{\widetilde{g}} \otimes
  \chi$ modulo $1+i$ satisfy the conditions for $\rho$. Consequently, when modulo $1+i$, the representations $
\rho^{\text{new}}_{+}$ and $\rho_{\widetilde{g}} \otimes
  \chi$ are surjective on $\GL_2(\F_2)$ and isomorphic up to semisimplification.

\medskip

We now find all $S_4$ extensions $E$ of $\Q(i)$ which contain $K$
and are unramified outside away from $(1+i)$ and $3$. There is a
unique order 48 group which has a normal subgroup isomorphic to
$S_4$, which is $S_4 \times \Z_2$ (The group $S_4 \rtimes \Z_2$ is
isomorphic to $S_4 \times \Z_2$ since all automorphisms of $S_4$ are
inner). Consequently, each $S_4$ extension of $\Q(i)$ corresponds to
a unique $S_4$ extension $L_4$ of $\Q$ which is unramified away from
2 and 3. Furthermore since $S_4\times \Z_2$ has a unique normal
subgroup isomorphic to $S_3$, $E$ has a unique subfield $L_3$  which
is Galois over $\Q$ with $\Gal(L_3/\Q)\cong S_3$. Over $\Q(i)$,
$f(x)$ can be obtained from $h(x)=x^3+3x-2$ (as
$(-x/i)^3+(-3/i)x-2=-i(x^3-3x-2i)$). The polynomial $h(x)$ appeared
in the proof of   \cite[Lemma 6.2]{lly05}.  Thus all such $S_4$
extensions $L_4$ of $\Q$ unramified away from 2 and 3 and containing
the splitting field of $h(x)$ are listed in \cite[pp. 138]{lly05}.
Accordingly, the following are the quartic polynomials which give
rise to $S_4$ extensions of $\Q(i)$ that contain $K$, unramified
outside of $(1+i)$ and $3$.
\begin{center}
\begin{tabular}{c c c}
defining equation     &  discriminant & $v$ with order $4$ Frobenius\\
\hline
$x^4-4ix-3=0$  &  $(3i)(1+i)^{18}(3)^2$ & $3+2i$ \\
$x^4-8ix+6=0$  &  $-(3i)(1+i)^{22}(3)^2$ & $3+2i$ \\
$x^4+12x^2-16ix+12=0$ & $(3i)(1+i)^{34}(3)^2$ & $6+i$\\
\end{tabular}
\end{center}

To conclude  the claim of Theorem \ref{thm:Ca}, it suffices to
compare the local Euler-$p$ factors of $L$-functions attached to the
two representations at $p=13,37,53$.

\subsection{The three-term Atkin and Swinnerton-Dyer
congruences}

We have remarked that in this case $ F\ta F^{-1}=\zeta ^{3(p-1)/2}
\ta$ and $F\zeta=\zeta^p F$. Let $$B_{-1}=\ta, \, B_{-3}=(\zeta
-\zeta^{-1}), \, B_3=\ta(\zeta -\zeta^{-1}).$$ They satisfy that
$B_{-1}, B_{-3}$, and $B_3$ are defined over $\Q_p$ when $p=5,7,11$
modulo 12 respectively. Moreover on $V^{\text{new}}$,
$$B_{-1}^2=-1, B_{-3}^2=B_3^2=-3.$$

Now apply the analysis in Section 7 of \cite{all05} with the roles
of $B_{\pm2}$ replaced by $B_{\pm 3}$ and $A$ replaced by $\ta$, we
obtain the following theorem.
\begin{theorem}[Atkin and Swinnerton-Dyer congruence for $S_3(\G_6)$]\label{Atkin and Swinnerton-Dyercong}
The Atkin and Swinnerton-Dyer congruence holds on the space
$S_3(\G_6) = \langle h_1,$ $ h_2, $ $ h_3, $ $ h_4, $ $h_5\rangle $.
More precisely, $h_3$ lies in $S_3(\G_2)$ and it satisfies the Atkin
and Swinnerton-Dyer congruence relations with the congruence form
$g_2(z)= \eta(4z)^6$; $h_2,h_4\in S_3(\G_3)$ and $h_2\pm i h_4$
satisfy Atkin and Swinnerton-Dyer congruence relations with level 27
congruence newforms (as in \cite{lly05}). For each odd prime $p\nmid
6l$, the subspace $\langle h_1, h_5\rangle $ has a basis depending
on the residue of $p$ modulo $12$ satisfying a 3-term Atkin and
Swinnerton-Dyer congruence at $p$ as follows.
\begin{enumerate}
\item If $p \equiv 1 \text{ mod }12$, then both $h_1$ and $h_3$ satisfy
the 3-term Atkin and Swinnerton-Dyer congruence at $p$ given by $T^2
- A(p)T + p^2$;
\item If $p \equiv 5 \text{ mod }12$, then $h_1$ (resp. $h_3$)  satisfies
the 3-term Atkin and Swinnerton-Dyer-congruence at $p$
 given by $T^2 - A(p)iT - p^2$ (resp. $T^2 + A(p)iT - p^2$);
 \item If $p \equiv 7 \text{ mod }12$ (resp. $p \equiv 11 \text{ mod }12$), then $h_1\pm h_3$
 (resp. $h_1\pm ih_3$) satisfy the 3-term
 Atkin and Swinnerton-Dyer congruence at $p$ given by
  $T^2 \mp \sqrt{-3}A(p)T -p^2$, respectively.
  \end{enumerate}
  All coefficients $A(p)\in \Z$ can be derived from the coefficients
  of $ \widetilde{g}$ via the modification described in Subsection
  \ref{subsec:g}.
  \end{theorem}

\section{Acknowledgements}
The author is  indebted to A.O.L. Atkin for his many original
observations and insightful comments. As a matter of fact, his
discoveries led to the discussions  in section \ref{sec:G6}.  The
author would also like to thank Wenching W. Li and Siu-Hung Ng for
their continuous valuable discussions and Chris Kurth for his
comments on an earlier version of this paper. Special thanks are due
to William A. Stein for computation assistance.


\newcommand{\etalchar}[1]{$^{#1}$}
\providecommand{\bysame}{\leavevmode\hbox
to3em{\hrulefill}\thinspace}
\providecommand{\MR}{\relax\ifhmode\unskip\space\fi MR }
\providecommand{\MRhref}[2]{%
  \href{http://www.ams.org/mathscinet-getitem?mr=#1}{#2}
} \providecommand{\href}[2]{#2}

\end{document}